\documentclass[a4paper,11pt]{article}
\usepackage{geometry}
\usepackage{mathtext}
\usepackage[latin2]{inputenc}
\usepackage[english]{babel}
\usepackage{amsmath}
\usepackage{amsfonts}
\usepackage{amssymb}
\usepackage{mathrsfs}
\usepackage{xcolor}
\usepackage{amsthm}
\usepackage{euscript}
\usepackage{authblk}

\DeclareMathOperator{\TT}{\mathbb{T}}
\DeclareMathOperator{\beps}{\boldsymbol{\varepsilon}}

\newcommand{\eLd}[1]{\left\|{#1}\right\|_{L_1(\mathbb{T}^d)}}

\renewcommand{\leq}{\leqslant}
\renewcommand{\geq}{\geqslant}

\theoremstyle{plain}\newtheorem{T1}{Theorem}
\theoremstyle{plain}
\theoremstyle{plain}
\theoremstyle{plain}
\theoremstyle{plain}
\theoremstyle{plain}\newtheorem{Cor}[T1]{Corollary}
\theoremstyle{plain}
\theoremstyle{definition}\newtheorem{Rem}[T1]{Remark}

\title{On Bernstein type quantitative estimates for Ornstein non-inequalities}
\author[1]{Krystian Kazaniecki}
\author[2]{Micha{\l} Wojciechowski}
\affil[1]{Institute of Analysis, Johannes Kepler University Linz}
\affil[1]{Institute of Mathematics, University of Warsaw}
\affil[2]{Institute of Mathematics, Polish Academy of Sciences}
    \date{}

\begin{document}

\maketitle
\begin{abstract} For the sequence of multi-indexes $\{\alpha_i\}_{i=1}^{m}$ and $\beta$ we study the inequality 
\[
\eLd{D^{\beta} f}\leq K_N \sum_{j= 1}^{m} \eLd{D^{\alpha_j}f},
\]
where $f$ is a trigonometric polynomial of degree at most $N$ on $d$-dimensional torus. Assuming some natural geometric property of the set $\{\alpha_j\}\cup\{\beta\}$ we show that
\[
K_{N}\geq C \left(\ln N\right)^{\phi},
\]
where $\phi<1$ depends only on the set $\{\alpha_j\}\cup\{\beta\}$.
\end{abstract}
\let\thefootnote\relax
\let\thefootnote\relax

\vskip 3mm
\section*{Introduction}In his inspiring article \cite{MR0149331} D. Ornstein showed that if $Q(D), P_1(D),\ldots  P_m(D)$ are homogeneous differential operators of the same order and $Q\notin \operatorname{span}\{P_j\} $ then, for any $C>0$, the inequality 
\begin{equation}\label{orn}
\eLd{Q(D) f}\leq C \sum_{j= 1}^{m} \eLd{P_j(D)f}
\end{equation}
 does not hold.
In particular, for any $C > 0$ and any multi-indices $\beta, \alpha_1, \ldots, \alpha_m$ with $|\beta|=|\alpha_1|=\ldots=|\alpha_m|$, $\beta\notin\{\alpha_j\}_{j=1}^{m}$
the inequality
\begin{equation}\label{ornmon}
\eLd{D^{\beta} f}\leq C \sum_{j= 1}^{m} \eLd{D^{\alpha_j}f}
\end{equation}
does not hold (in this paper $L_1(\TT^d)$ is considered with respect to the normalized Haar measure).

In this paper we deal with the quantitative version of this theorem. We are interested in the constant of Bernstein type i.e.  what is the growth of the best constant C in \eqref{orn} when the inequality is restricted to the polynomials of degree at most $n$. To the best of our knowledge no result of such type is known. 

Our  main idea is to use the properties of finite Riesz products \cite{MR1544321}.
In fact we are constructing explicitly the trigonometric polynomials for which our bounds hold.  For different (but qualitative) proofs of Ornstein non-inequality in the isotropic case check for example \cite{C4},\cite{KK6}.

Focusing for a moment on the simplest case of our results we get the following:

\begin{Cor}\label{prosty}
For every $N\in\mathbb{N}$ there exists a trigonometric polynomial $P_{N}$ on $\mathbb{T}^2$ of degree $N$, which satisfies  
\[
\eLd{\frac{\partial^2}{\partial x^2} P_{N}}+\eLd{\frac{\partial^2}{\partial y^2} P_{N}}\leq 1
\]
but 
\[
\eLd{\frac{\partial^2}{\partial x\partial y} P_{N}} > C \ln^{\frac{1}{2}} N.
\]
\end{Cor}

We do not know if the bound from Corollary \ref{prosty} is sharp. In fact one can establish in a rather standard way that the mixed derivative from Corollary 1 could not have norm greater than $c\ln n$. Indeed the (linear and invariant) operator $T$ which retrieves the mixed derivative from the pure ones is of a weak type $(1,1)$ (see \cite{Fabes1966}). Hence, by the Nikolskii type inequality for Lorentz spaces (see \cite[Theorem 3]{She84} and \cite[Lemma 3.1]{Aki}),  for a trigonometric polynomial $f$ of degree $N$
\[
\begin{split}
\|\frac{\partial^2}{\partial x \partial y} f\|_{L_{1}(\mathbb{T}^2)}&=\|Tf\|_{L_{1,1}(\mathbb{T}^2)}\leq C\ln(1+ N) \|Tf\|_{L_{1,\infty}(\mathbb{T}^2)}\\ &\leq C\ln(1+ N)\left(\|\frac{\partial^2}{\partial x^2} f\|_{L_{1}(\mathbb{T}^2)}+\|\frac{\partial^2}{\partial y^2} f\|_{L_{1}(\mathbb{T}^2)}\right).
\end{split}
\] 
The same comment concerns all the results obtained in this paper. All bounds from below presented here are of the form $(\ln N)^{\phi}$ for some $\phi<1$, while the common upper bound is $\ln N$. Nevertheless we conjecture that the optimal exponent $\phi$ should be equal to one (see Remark \ref{ostatni}).

This paper contains final results of the study we began in \cite{KW}. It seems it is the first use in the literature of trigonometric polynomials in the context of Ornstein non-inequalities.

\subsection*{Acknowledgments} Authors would like to thank Maciej Rzeszut for helpful discussions. We would like to express our gratitude to anonymous reviewers whose comments/suggestions helped improve and clarify this manuscript.

This research was supported by the National Science Centre, Poland, and Austrian Science Foundation FWF joint CEUS programme. FWF project no. I5231 and National Science Centre project no. 2020/02/Y/ST1/00072.

\section*{Results}
In this paper we consider more general, anisotropic case of the inequality \eqref{ornmon} i.e. such that there exists $\Lambda \in  \mathbb{N}^d$ with $\langle\alpha_1\;,\;\Lambda\rangle=\langle\alpha_2\;,\;\Lambda\rangle=\cdots=\langle\beta\;,\;\Lambda\rangle$. This case was already considered in the literature (see \cite{KSW}), with an additional assumption on the parity of derivatives $|\alpha_j|\equiv_2|\beta|$. (We mention that the results obtained there were only of the qualitative nature.) In our present approach we remove this "parity assumption". However we still need other geometric conditions.

The paper contains two results which are proved in quite similar way. Each of them provides a geometric criterion for a set of symbols of partial derivatives which yields quantitative estimates of Ornstein type.

\begin{T1}\label{glownetw1}
 Assume $\alpha_1, \ldots, \alpha_n,\beta$ are multi-indices in  $\left(\mathbb{N}\cup\{0\}\right)^d$ and there exists a pair of vectors $\Gamma,\; \Lambda\in\mathbb{N}^d$ for which the following conditions are satisfied
 \[
 \langle \alpha_1,  \Lambda \rangle= \langle \beta,  \Lambda \rangle=\langle \alpha_2,  \Lambda \rangle=\ldots= \langle \alpha_m,  \Lambda \rangle,
 \]
and
\begin{equation}\label{gammaporz}
\langle \alpha_1,  \Gamma \rangle > \langle \beta,  \Gamma \rangle>\langle \alpha_2,  \Gamma \rangle\geq\ldots\geq\langle \alpha_m,  \Gamma \rangle,
\end{equation}
Let $K_N$  be the smallest constant such that for every trigonometric polynomial $f$ of degree at most $N$ the following estimate hold 
\begin{equation}\label{glowna1}
\eLd{D^{\beta} f}\leq K_N \sum_{j= 1}^{m} \eLd{D^{\alpha_j}f}.
\end{equation}
Then, there exists a constant $C>0$ such that 
\[
K_N>C (\ln N)^{\phi},
\]
where $\phi=\frac{1}{2}\left(1 - \frac{\langle\alpha_1-\beta ,  \Gamma \rangle}{\langle\alpha_1-\alpha_2 ,  \Gamma \rangle}\right)$.
\end{T1}
\begin{Rem}
The inequalities \eqref{gammaporz} could be satisfied for different vectors $\Gamma$ and then for different permutations of the set $\{\alpha_1,\ldots,\alpha_m\}$. For any fixed set of multi-indices the choice of the optimal vector $\Gamma$ is a simple  optimization problem. In dimension 2 for fixed $\alpha_1$ the value of $\phi$ doesn't depend on choice of $\Gamma$.  
\end{Rem}
\begin{T1}\label{glownetw2}
Assume $\alpha_1, \ldots, \alpha_n,\beta$ are multi-indices in  $\left(\mathbb{N}\cup\{0\}\right)^d$ and there exists a vector $\Lambda\in\mathbb{N}^d$\footnote{In this paper we put $\mathbb{N}=\{1,2,3,\ldots\}$} for which the following condition is satisfied
 \[
 \langle \alpha_1,  \Lambda \rangle= \langle \beta,  \Lambda \rangle=\langle \alpha_2,  \Lambda \rangle=\ldots= \langle \alpha_m,  \Lambda \rangle.
 \]
Suppose moreover that there exists $\beps\in\{0,1\}^d$ such that
 \begin{equation}\label{epscal}
   \langle \beta, \beps  \rangle \not\equiv  \langle \alpha_1, \beps  \rangle\; \mbox{mod } 2
  \qquad \mbox{and}\qquad
    \langle \alpha_{j}, \beps \rangle \equiv  \langle \alpha_1, \beps  \rangle\; \mbox{mod } 2 
 \end{equation}
 for $j\in\{1, \ldots , m\}$.
 Let $K_N$  be the smallest constant such that, for every trigonometric polynomial $f$ of degree at most $N$, the following estimate hold 
\begin{equation}\label{glowna2}
\eLd{D^{\beta} f}\leq K_N \sum_{j= 1}^{m} \eLd{D^{\alpha_j}f}.
\end{equation}
Then, there exists a constant $C>0$ such that 
\[
K_N>C (\ln N)^{\frac{1}{2}}.
\]
\end{T1}
\begin{Rem}
The case $\beps=(1,1,\ldots,1)$ corresponds to the anisotropic Sobolev space which contains an invariant, complemented, infinite dimensional subspace isomorphic to  a Hilbert space (see \cite{Pelczynski1992} for details).
\end{Rem}

\section*{Proof of Theorem \ref{glownetw1}}
\begin{proof}
Let $\Lambda=(\lambda_1,\ldots,\lambda_d)$ and $\Gamma=(\gamma_1,\ldots,\gamma_d)$.  We introduce an auxiliary sequence $(b_k)_{k\geq 1}$ depending on the parity of our multi-indices. If $|\alpha_1| -|\beta|$ is even we put $b_k:= (2+(-1)^k)$, if not we put $b_k:=1$.  For a fixed $n\in\mathbb{N}$ we define a sequence of vectors $(a_k)_{k\geq 1}$ in $\mathbb{N}^d$ by the formula $a_k=( a_k(1),\ldots, a_k(d))$, where
\begin{equation}\label{defin}
a_k(j):= 3^{\lambda_j 2 k n} b_k^{\gamma_j} \lfloor  n^{\theta \gamma_j }\rfloor
\end{equation}
and
\begin{equation}\label{deftheta}
\theta:=\frac{1}{\langle\alpha_1-\alpha_2 , \Gamma\rangle} 
\end{equation} 
Since $\gamma_j$ and $\lambda_j$ are positive integers we know that, for any $j\in\{1,2,\ldots, d\},$ 
\begin{equation}\label{nawspol}
a_k(j)>3^{2(n-1)}a_{k-1}(j).
\end{equation}
and
\begin{equation}\label{przyrost}
\|a_k\|_2> 3^{2(n-1)} \|a_{k-1}\|_2.
\end{equation} 
We define a modified Riesz product based on this sequence
 \begin{equation}\label{Riesze}
R_n(x)=-1 + \prod_{k=1}^{n}\left(1+\cos\langle x  ,  a_k\rangle \right),
 \end{equation}
and the family of sets 
\begin{equation*}
A_k=\{q : q=  a_k + \sum_{j=1}^{k-1} \xi_j a_j,\; \xi_j\in\{-1,0,1\} \mbox{ for } j\in\{1,\ldots, k-1\}\;\}. 
\end{equation*}
From \eqref{przyrost} by standard calculations we know that every point in $A_k$ has a unique representation as $ a_k + \sum\limits_{j=1}^{k-1} \xi_j a_j$.  From \eqref{nawspol} there exists a constant $2\geq\tau>1$ independent of $k$ and $j$ such that 
\begin{equation}\label{oszacowaniewspol}
\frac{1}{\tau} a_k(j)\leq |q(j)|\leq \tau a_k(j), 
\end{equation}
for all $q\in A_k$.

For $\mu \in \mathbb{Z}^d$ we denote $n^{\mu}=\prod_{j=1}^{d} n_j^{\mu(j)}$. For $q\in A_k$ of the form $q= a_k+\sum_{j=1}^{k-1} \xi_j a_j$ we set $r(q)=\# \{j: \xi_j\neq 0\}+1$ and $r(-q)=r(q)$. Let $W_n(x)$ be a polynomial given by the formula
\begin{equation}\label{defW}
W_n(x)= \sum_{k=1}^{n}\sum_{q\in A_k\cup -A_k} \frac{i^{-|\alpha_1|}}{q^{\alpha_1}}\frac{1}{2^{r(q)}} e^{i\langle q,  x\rangle}.
\end{equation}
Note that
\[
D^{\alpha_1} W_n(x)= \sum_{k=1}^{n}\sum_{q\in A_k\cup -A_k} \frac{1}{2^{r(q)}} e^{i\langle q,  x\rangle}=R_n(x).
\]
Moreover, for $\mu\in\{\beta, \alpha_2,\ldots, \alpha_m\}$
\[
D^{\mu} W_n(x)=\sum_{k=1}^{n}\sum_{q\in A_k\cup -A_k}i^{|\mu|-|\alpha|}\frac{q^{\mu}}{q^{\alpha_1}} \frac{1}{2^{r(q)}} e^{i\langle q,  x\rangle},
\]
which could be represented as
\begin{equation}\label{rozklad}
D^{\mu} W_n(x)= B_{\mu,n}(x) + G_{\mu,n}(x),
\end{equation}
where
\begin{equation}\label{defB}
  B_{\mu,n}(x)=\sum_{k=1}^{n}\sum_{q\in A_k}\frac{i^{|\mu|-|\alpha_1|}}{2^{r(q)}} \left(\left(\frac{q^{\mu}}{q^{\alpha_1}}-\frac{a_{k}^{\mu}}{a_k^{\alpha_1}}\right) e^{i\langle q,  x\rangle}+\left(\frac{(-q)^{\mu}}{(-q)^{\alpha_1}}-\frac{(-a_{k})^{\mu}}{(-a_k)^{\alpha_1}}\right)  e^{i\langle -q,  x\rangle}\right)  
\end{equation}
and
\begin{equation}\label{defG}
G_{\mu,n}(x) =\sum_{k=1}^{n}i^{|\mu|-|\alpha_1|}\frac{a_{k}^{\mu}}{a_k^{\alpha_1}}\sum_{q\in A_k}\frac{1}{2^{r(q)}}\left( e^{i\langle q,  x\rangle}+ (-1)^{|\mu|-|\alpha_1|}e^{i\langle -q,  x\rangle}\right).
\end{equation}
First we estimate the $L_1$-norm of $B_{\mu,n}$. Let $v=(\frac{q(1)}{a_k(1)},\ldots,\frac{q(d)}{a_k(d)})$ for $q\in A_k$. From \eqref{nawspol}
\begin{equation}\label{raznizej}
\|v-\mathbf{1}\|_{2}\leq C(d) 3^{-2n},
\end{equation}
where $\mathbf{1}:=(1,\dots,1)$.
Observe that for $q\in A_k$ by \eqref{oszacowaniewspol}, \eqref{raznizej} and by Lipschitz continuity of functions $x^{\alpha_1}$, $x^{\mu}$ on the cube $[\frac{1}{\tau}, \tau]^{d}$ we get
\begin{equation}\label{eq: lipest}
 \begin{split}
 \big|q^{\mu}a_k^{\alpha_1}-a_k^{\mu}q^{\alpha_1}\big|&\leq \big|q^{\mu}\left(a_k^{\alpha_1}-q^{\alpha_1}\right)\big|+ \big|q^{\alpha_1}\left(q^{\mu}- a_k^{\mu}\right)\big|
 \\&\leq C \left(|q^{\alpha_1}||a_k^{\mu}|\big| 1-v^{\alpha_1}\big|+|q^{\mu}||a_k^{\alpha_1}|\big|1-v^{\mu}\big|\right)
 \\&\stackrel{\eqref{oszacowaniewspol}}{\leq} C |a_k^{\alpha_1}||a_k^{\mu}|\left( \big| 1-v^{\alpha_1}\big|+\big|1-v^{\mu}\big|\right) 
 \\&\stackrel{Lip.}{\leq} C \|\mathbf{1}-v\|_{2} |a_k^{\alpha_1}||a_k^{\mu}|
 \\&\stackrel{\eqref{raznizej}}{\leq} C 3^{-2n} |a_k^{\alpha_1}||a_k^{\mu}|.
 \end{split}    
\end{equation}
We calculate $a_k^{\mu}$ and $a_k^{\alpha_1}:$
\begin{equation}\label{eee}
      a_k^{\mu}=3^{\sum\limits_{j=1}^{d} \mu(j) \lambda_j 2k n} b_k^{\sum\limits_{j=1}^{d} \mu(j) \gamma_j } \prod_{i=1}^{d} \lfloor n^{\theta\gamma_i}\rfloor^{\mu(j)}= 3^{\langle\mu, \Lambda\rangle2 kn} b_k^{\langle\mu, \Gamma\rangle}\prod_{j=1}^{d} \lfloor n^{\theta\gamma_j}\rfloor^{\mu(j)},
\end{equation}
and similarly (replacing $\mu$ by $\alpha_1$),
\begin{equation}\label{eee1}
a_k^{\alpha_1}= 3^{\langle\alpha_1, \Lambda\rangle2 kn}b_k^{\langle\alpha_1, \Gamma\rangle} \prod_{j=1}^{d} \lfloor n^{\theta\gamma_j}\rfloor^{\alpha_1(j)}.      
\end{equation}

Since we only use a finite number of exponents, there is a constant $C>1$ such that for any $\nu\in\{\beta,\alpha_1,\ldots,\alpha_m\},$
\begin{equation} \label{podlg}
\frac{1}{C}n^{\theta\cdot\langle \nu , \Gamma\rangle} \leq\prod_{j=1}^{d}\lfloor n^{\theta\gamma_j}\rfloor^{\nu(j)}\leq Cn^{\theta\cdot\langle \nu , \Gamma\rangle}.
\end{equation}
For every $\mu$ as above we have $\langle\mu , \Lambda\rangle= \langle\alpha_1 , \Lambda\rangle$. By \eqref{oszacowaniewspol}, \eqref{eq: lipest} and \eqref{eee}, \eqref{eee1}, \eqref{podlg} we get
\[
\big|\left(\frac{q^{\mu}}{q^{\alpha_1}}-\frac{a_{k}^{\mu}}{a_k^{\alpha_1}}\right)\big|=\big|\frac{q^{\mu}a_k^{\alpha_1}-a_k^{\mu}q^{\alpha_1}}{q^{\alpha_1} a_k^{\alpha_1}}\big|\leq C 3^{-2n}\frac{|a_k^{\alpha_1}||a_k^{\mu}|}{|a_k^{\alpha_1}|^2}\leq C 3^{-2n} n^{\langle\mu-\alpha_1, \Gamma\rangle}\leq C 3^{-n}.
\]
 Plugging the above estimates for $\mu\in\{\beta, \alpha_2,\ldots, \alpha_m\}$ into the formula for $B_{\mu,n}$ we get 
\[
\eLd{B_{\mu,n}}\leq \sum_{k=1}^{n} \sum_{q\in A_k} 2C 3^{-n}\leq 2C 3^{-n} 3^n=2C.
\]

We pass to the estimates of the $L_1$ norm of $G_{\alpha_j,n}$ for $j\geq 2$. For $k\in \mathbb{N}$ and $1\leq k\leq n$ we define 
\[
\psi_{k}(x)=\prod_{l=1}^{k-1}\left(1+\cos\langle x ,    a_{l}\rangle \right),
\]
Simple algebraic manipulation gives us
\[
 G_{\alpha_j,n}(x)=\sum_{k=1}^{n}i^{|\mu|-|\alpha_1|}\frac{a_{k}^{\alpha_j}}{a_k^{\alpha_1}}\frac{1}{2}\left(e^{i\langle a_k,   x\rangle}+ (-1)^{|\alpha_j|-|\alpha_1|}e^{i\langle -a_k,   x\rangle}\right)\psi_k(x).
\]
Since $\langle \alpha_j ,  \Lambda\rangle= \langle \alpha_1 ,  \Lambda\rangle$, by \eqref{eee} and \eqref{eee1} we get   
\[
\bigg|\frac{a_{k}^{\alpha_j}}{a_k^{\alpha_1}}\bigg|\leq C n^{\theta\cdot\langle \alpha_j-\alpha_1 ,  \Gamma\rangle}.
\]
Therefore 
\[
\eLd{G_{\alpha_j,n}}\leq C n^{\theta\cdot\langle \alpha_j-\alpha_1,  \Gamma\rangle} \sum_{k=1}^{n} \eLd{\psi_k}.
\]
As Riesz products $\psi_l$'s have  $L_1$ norms equal to 1. Hence
\[
\eLd{G_{\alpha_j,n}}\leq C n^{\theta\cdot\langle \alpha_j-\alpha_1,  \Gamma\rangle +1 }.
\]
By \eqref{gammaporz} and \eqref{deftheta}, for any $j\in\{2,\ldots,m\},$ 
\[
\theta\cdot \langle \alpha_j-\alpha_1 ,  \Gamma\rangle =\frac{\langle \alpha_j-\alpha_1 ,  \Gamma\rangle}{\langle \alpha_1-\alpha_2 ,  \Gamma\rangle}\leq - 1.
\]
Hence there exists $C>0$ such that
\[
\eLd{G_{\alpha_j,n}}\leq  C,
\]
for any $n\in \mathbb{N}$ and any $j\in\{2,\cdots,m\}$.
Therefore for $j\in\{2,\cdots, m\}$ and $n\in\mathbb{N}$ 
\[
\eLd{D^{\alpha_j} W_n}\leq \eLd{B_{\alpha_j,n} }+\eLd{G_{\alpha_j,n} }\leq C.
\] 
Since $D^{\alpha_1} W$ is a modified Riesz product, 
\[
\eLd{D^{\alpha_1} W_n}\leq 2.
\]
Summing the above inequalities we get   
\begin{equation}\label{Prawast}
\sum_{j= 1}^{m} \eLd{D^{\alpha_j}W_n} \leq C.
\end{equation}
Now we estimate $\eLd{D^{\beta} W_n}$ from below. Since the norm of $B_{\beta,n}$ is uniformly bounded with respect to $n$, it is enough to show that the norm of $G_{\beta,n}$ is large. 

\begin{Rem}\label{Meyer}
 In \cite[Remark on p. 563]{MR0240563} Y. Meyer observes that the condition $\sum_{k=1}^{\infty} \frac{a_{k}(j)}{a_{k+1}(j)} < +\infty$ yields
\[
\begin{split}
\bigg\|\sum_{\xi\in \{-1,0,1\}^n} b\left(\sum_{k=1}^{n} \xi_{k} a_{k}(j)\right) \operatorname{exp}&\left(i \sum_{k=1}^{n} \xi_{k} a_{k}(j)t\right)   \bigg\|_{L^1(\mathbb{T})}
\\&\simeq \bigg\|\sum_{\xi\in \{-1,0,1\}^n}b\left(\sum_{k=1}^{n} \xi_{k} a_{k}(j)\right) \operatorname{exp}\left(i \sum_{k=1}^{n} \xi_{k} t_{k}\right)   \bigg\|_{L^1(\mathbb{T}^{n})}
\end{split}
\]
The constant in the above isomorphism depends only on the value of $\sum_{k=1}^{\infty} \frac{a_{k}(j)}{a_{k+1}(j)}$. For elementary proofs of this fact see note \cite{KWEQ} or \cite[Proposition 4]{Bonami2020}.
By a simple tensoring argument 
\[
\begin{split}
\bigg\|\sum_{\xi\in \{-1,0,1\}^n} b\left(\sum_{k=1}^{n} \xi_{k} a_{k}\right) \operatorname{exp}&\left(i \sum_{k=1}^{n}\langle \xi_{k} a_{k}  ,  t\rangle\right)   \bigg\|_{L^1(\mathbb{T}^d)}
\\&\simeq \bigg\|\sum_{\xi\in \{-1,0,1\}^n}  b\left(\sum_{k=1}^{n} \xi_{k} a_{k}\right) \operatorname{exp}\left(i \sum_{k=1}^{n} \langle\xi_{k}  ,   t_{k}\rangle\right)   \bigg\|_{L^1(\mathbb{T}^{{nd}})}
\end{split}
\]
In our case  there exists a constant $C>0$  independent of $n$ such that finite sequence $(a_k(j))_{k=1}^{n}$ define by \eqref{defin} satisfies
\[
\sum_{k=1}^{n} \frac{a_{k}(j)}{a_{k+1}(j)} \stackrel{\eqref{nawspol}}{<} n 3^{-2n}< C< +\infty
\] 
for any $j\in\{1,\ldots,d\}$. Hence when calculating  $L_1$  norm we can treat exponents with different $a_k$ as independent random variables.  
\end{Rem}
We consider two cases separately.
\subsection*{Case I: $|\alpha_1|-|\beta|$ is even.} In this case 
\[
|G_{\beta,n}(x)|= \bigg|\sum_{k=1}^{n}\frac{a_{k}^{\beta}}{a_k^{\alpha_1}}\cos(\langle a_k, x\rangle)\psi_k(x)\bigg|=\bigg|\frac{a_{1}^{\beta}}{a_1^{\alpha_1}}\psi_{1}+ \frac{a_{n}^{\beta}}{a_n^{\alpha_1}}\psi_{n+1}+ \sum_{k=1}^{n-1} \left(\frac{a_{k+1}^{\beta}}{a_{k+1}^{\alpha_1}}-\frac{a_{k}^{\beta}}{a_k^{\alpha_1}}\right) \psi_{k+1}\bigg|. 
\]
Since $\psi_l$ are Riesz products, by \eqref{przyrost} and by the inequality of Lata{\l}a (Theorem 1 \cite{Latala2014}, see also \cite{Bonami2020}, \cite{Damek2015}).
\[
\eLd{G_{\beta,n}}\geq \bigg|\frac{a_{1}^{\beta}}{a_1^{\alpha_1}}\bigg|+ \bigg|\frac{a_{n}^{\beta}}{a_n^{\alpha_1}}\bigg|+ \sum_{k=1}^{n-1}\bigg| \frac{a_{k+1}^{\beta}}{a_{k+1}^{\alpha_1}}-\frac{a_{k}^{\beta}}{a_k^{\alpha_1}} \bigg|. 
\]
From the definition of $a_k$ (and the fact that $\langle \alpha_1,  \Lambda \rangle= \langle \beta,  \Lambda \rangle$) we get
\[
\bigg| \frac{a_{k+1}^{\beta}}{a_{k+1}^{\alpha_1}}-\frac{a_{k}^{\beta}}{a_k^{\alpha_1}} \bigg| \geq C n^{\theta\cdot\langle \beta-\alpha_1 ,  \Gamma\rangle } |b_k^{\langle \beta ,  \Gamma\rangle}-b^{\langle \beta ,  \Gamma\rangle}_{k+1}|\geq Cn^{\theta\cdot\langle \beta-\alpha_1 ,  \Gamma\rangle}|3^{\langle \beta ,  \Gamma\rangle}-1|.
\]
From the definition of $\theta$ we get
\[
\eLd{G_{\beta,n}}\geq C n^{1 - \frac{\langle\beta-\alpha_1 ,  \Gamma \rangle}{\langle\alpha_2-\alpha_1 ,  \Gamma \rangle}}
\]
and consequently
\[
\eLd{D^{\beta} W_n}\geq C n^{1 - \frac{\langle\beta-\alpha_1 ,  \Gamma \rangle}{\langle\alpha_2-\alpha_1 ,  \Gamma \rangle}}.
\]
There exists $C>0$ independent of $n$ such that $|a_k|\leq 3^{Cn^2}$ for $1\leq k\leq n$. Therefore $\operatorname{deg} W_n(x)\leq 3^{C n^2}$ and $\ln(\operatorname{deg} W_n(x))\leq C n^2$. Hence for large enough $n$ there exists a constant $C>0$ such that 
\[
\eLd{D^{\beta} W_n}\geq C (\ln \operatorname{deg} W_n )^{\frac{1}{2}\left(1 - \frac{\langle\beta-\alpha_1 ,  \Gamma \rangle}{\langle\alpha_2-\alpha_1 ,  \Gamma \rangle}\right)}.
\]
From \eqref{Prawast} and \eqref{glowna1} we get 
\[
K\geq C (\ln \operatorname{deg} W_n )^{\frac{1}{2}\left(1 - \frac{\langle\beta-\alpha_1 ,  \Gamma \rangle}{\langle\alpha_2-\alpha_1 ,  \Gamma \rangle}\right) }.
\]
\subsection*{Case II: $|\alpha_1|-|\beta|$ is odd.} 
In this case 
\[
|G_{\beta,n}(x)|= \bigg|\sum_{k=1}^{n}\frac{a_{k}^{\beta}}{a_k^{\alpha_1}}\sin(\langle a_k, x\rangle)\psi_k(x)\bigg|.
\]
By the definition of the sequence $a_k$ 
\[
\left|\frac{a_{k}^{\beta}}{a_k^{\alpha_1}} - n^{ - \frac{\langle\beta-\alpha_1 ,  \Gamma \rangle}{\langle\alpha_2-\alpha_1 ,  \Gamma \rangle}}\right| = \left|\frac{\prod_{s=1}^{d} \lfloor n^{\theta\gamma_j}\rfloor^{\beta(s)}}{\prod_{s=1}^{d} \lfloor n^{\theta\gamma_j}\rfloor^{\alpha_1(s)}} - n^{ - \frac{\langle\beta-\alpha_1 ,  \Gamma \rangle}{\langle\alpha_2-\alpha_1 ,  \Gamma \rangle}}\right|\leq \frac{C}{n^{\theta}}\cdot n^{ - \frac{\langle\beta-\alpha_1 ,  \Gamma \rangle}{\langle\alpha_2-\alpha_1 ,  \Gamma \rangle}}.
\]
Therefore
\[
\eLd{G_{\beta,n}}\geq n^{ - \frac{\langle\beta-\alpha_1 ,  \Gamma \rangle}{\langle\alpha_2-\alpha_1 ,  \Gamma \rangle}}\eLd{\sum_{k=1}^{n}\sin(\langle a_k, x\rangle)\psi_k(x)} -\frac{C}{n^{\theta}} \cdot n^{1 - \frac{\langle\beta-\alpha_1 ,  \Gamma \rangle}{\langle\alpha_2-\alpha_1 ,  \Gamma \rangle}}.
\]
However
\[
\begin{split}
\eLd{\sum_{k=1}^{n}\sin(\langle a_k, x\rangle)\psi_k(x)}&=\eLd{\sum_{k=1}^{n}\left(\cos(\langle a_k, x\rangle)-e^{i\langle a_k, x\rangle}\right)\psi_k(x)}
\\&\geq \eLd{\sum_{k=1}^{n}e^{i\langle a_k, x\rangle}\psi_k(x)} - \eLd{\psi_{n+1}}
\\&\geq \eLd{\sum_{k=1}^{n}e^{i\langle a_k, x\rangle}\psi_k(x)} -1.
\end{split}
\]
The sequence $a_k$ satisfies the assumptions of Meyer's theorem (see Remark \ref{Meyer}). Because of that we can use Lemma 2 from \cite{MR1649869}, which gives 
\begin{equation*}
\left\|\sum_{j=1}^{n} e^{ i \langle a_k,x\rangle}\psi_k(x)\right\|_{L_1(\mathbb{T}^d)}\geq Cn.
\end{equation*} 
Thus
\[
\eLd{G_{\beta,n}}\geq \left(C_1 -\frac{C}{n^{\theta}}\right)n^{1 - \frac{\langle\beta-\alpha_1 ,  \Gamma \rangle}{\langle\alpha_2-\alpha_1 ,  \Gamma \rangle}}.
\]
Therefore for large enough $n$ we get 
\[
\eLd{G_{\beta,n}}\geq Cn^{1 - \frac{\langle\beta-\alpha_1 ,  \Gamma \rangle}{\langle\alpha_2-\alpha_1 ,  \Gamma \rangle}}
\]
and similarly as is in Case I we obtain 
\[
K\geq C (\ln \operatorname{deg} W_n)^{\frac{1}{2}\left(1 - \frac{\langle\beta-\alpha_1 ,  \Gamma \rangle}{\langle\alpha_2-\alpha_1 ,  \Gamma \rangle}\right)} .
\]
\end{proof}
\section*{Proof of Theorem \ref{glownetw2}}
\begin{proof}
We prove Theorem \ref{glownetw2} in an analogous way. First we define the sequence $a_k$ by the formula
\[
a_k(j)= 3^{\lambda_j 2 k n} (-1)^{\beps_j k}.
\]
Once again we use the modified Riesz products \eqref{Riesze} and the corresponding polynomials $W_n$ \eqref{defW}. As in the proof of Theorem \ref{glownetw1} we have 
\[
D^{\mu} W_n(x)= B_{\mu,n}(x) + G_{\mu,n}(x)
\]
for any $\mu\in\{\beta,\alpha_2,\ldots ,\alpha_m\}$ and $B_{\mu,n}(x)$, $G_{\mu,n}(x)$ defined as in \eqref{defB} and \eqref{defG}. Since the sequence $a_k$ has super-exponential growth
\[
\frac{|a_k(j)|}{|a_{k+1}(j)|}\leq 3^{-2n},
\]
we obtain the bounds on $B_{\mu,n}(x)$
\[
\eLd{B_{\mu,n}}\leq C
\]
for any $\mu\in \{\beta,\alpha_2,\ldots ,\alpha_m\}$.
Note that by \eqref{epscal}
\[
a_k^{\alpha_j}=3^{\langle \Lambda,\; \alpha_j\rangle} (-1)^{k \langle\beps,\; \alpha_j\rangle }=3^{\langle \Lambda\;,\; \alpha_1\rangle} (-1)^{k \langle\beps,\; \alpha_1\rangle }= a_k^{\alpha_1}.
\]
Hence by \eqref{defG}, 
\[
G_{\alpha_j,n}(x) =\sum_{k=1}^{n}i^{|\alpha_j|-|\alpha_1|}\sum_{q\in A_k}\frac{1}{2^{r(q)}}\left( e^{i\langle q,  x\rangle}+ (-1)^{|\alpha_j|-|\alpha_1|}e^{i\langle -q,  x\rangle}\right).
\]
Thus for $|\alpha_1|\equiv |\alpha_j| \mod 2$ we get
\[
\eLd{G_{\alpha_j,n}}=\eLd{-1 + \prod_{k=1}^{n}\left(1+\cos\langle x  ,  a_k\rangle\right)}=\eLd{R_n}\leq 2
\]
and for $|\alpha_1|\equiv |\alpha_j|+1 \mod 2$
\[
\eLd{G_{\alpha_j,n}}=\eLd{-1 + \prod_{k=1}^{n}\left(1+\sin\langle x  ,  a_k\rangle\right)}\leq 2.
\]
The only thing left to do is the estimate on the norm of $G_{\beta,n}$ from below. By \eqref{epscal} we get 
\[
a_k^{\beta}=3^{\langle \Lambda,\; \beta\rangle} (-1)^{k \langle\beps,\; \beta\rangle }=3^{\langle \Lambda,\; \alpha_1\rangle} (-1)^{k (\langle\beps,\; \alpha_1\rangle+1) }= (-1)^{k} a_k^{\alpha_1}.
\]
Therefore 
\[
G_{\beta,n}(x) =\sum_{k=1}^{n} (-1)^{k} i^{|\beta|-|\alpha_1|}\sum_{q\in A_k}\frac{1}{2^{r(q)}}\left( e^{i\langle q,  x\rangle}+ (-1)^{|\beta|-|\alpha_1|}e^{i\langle -q,  x\rangle}\right).
\]
Let $g$ be given by the formula
\[
g_k(x)=\left\{\begin{array}{cr}
     \prod_{k=1}^{n}\left(1+\cos\langle x  ,  a_k\rangle\right),&\;\mbox{ for } |\beta|\equiv |\alpha_1|\mod 2, \\
    \prod_{k=1}^{n}\left(1+\sin\langle x  , a_k\rangle\right),  & \; \mbox{ for }  |\beta|\not\equiv |\alpha_1|\mod 2.
\end{array}
\right.
\]
Then
\[
|G_{\beta,n}(x)|= \left|\sum_{k=0}^{n-1} (-1)^k\left( g_{k+1}(x) -g_{k}(x)\right)\right|=\left|(-1)^{n-1}  g_n(x)+g_0(x)+ \sum_{k=1}^{n-1} 2 (-1)^{k-1} g_{k}(x)\right|.
\]
Applying  Lata{\l}a's inequality (Theorem 1 in \cite{Latala2014}) we obtain 
\begin{equation}\label{odwolanieMichal}
\eLd{G_{\beta,n}}\geq C n.    
\end{equation}
As in previous section we get 
\[
K\geq C (\ln \operatorname{deg} W_n)^{\frac{1}{2}}. 
\]
\end{proof}
\begin{Rem}
Actually to obtain estimate \eqref{odwolanieMichal} we could use a weaker (random) form of Lata{\l}a's inequality (see Lemma 1 in \cite{MR1649869}). However, to do this one needs to adjust the construction to the randomness of choice of signs, which significantly complicates the redaction. 
\end{Rem}
\begin{Rem}\label{ostatni}
The careful study of the above proofs shows that the source of sub-logarithmic growth of constant lies in super-exponential growth of the sequence $a_k$ (see \eqref{nawspol}). There are two reasons for significant growth of this sequence. In order to use the Lata{\l}a's inequality just geometrical growth would be enough (see \cite{Bonami2020}). However, we use a Riesz product as one of the derivatives involved in the proof. To recover the Riesz product structure for the remaining derivatives - which we need to use in the inequality of Lata{\l}a - we perturbed the actual functions. Our method to control the arising error terms requires super-exponential growth of the sequence $a_k$. It seems that any improvement of this method would require a more delicate study of the aforementioned error terms.

\end{Rem}
\bibliographystyle{siam}
\bibliography{cont}
Krystian Kazaniecki\\
Institute of Analysis, Johannes Kepler University Linz\\
Institute of Mathematics, University of Warsaw\\
krystian.kazaniecki@jku.at\\
 \\
Micha{\l} Wojciechowski\\
Institute of Mathematics, Polish Academy of Sciences\\
miwoj@impan.pl

\end{document}